\documentclass[a4paper]{article} \def\zibreport{1}

\usepackage{ifthen}

\usepackage{amsmath,amsfonts,amssymb}
\usepackage{xspace}
\usepackage{wrapfig}
\usepackage{url}
\usepackage{todonotes}
\usepackage[hyperfootnotes=false]{hyperref}

\ifthenelse{\zibreport = 1}{
  \let\pdfoutorg\pdfoutput
  \let\pdfoutput\undefined
  \usepackage{zibtitlepage}
  \let\pdfoutput\pdfoutorg

  \ZTPTitle{Exploring the Numerics of Branch-and-Cut for Mixed Integer Linear Optimization}
  \ZTPAuthor{Matthias~Miltenberger$^1$, Ted~Ralphs$^2$, Daniel~E.~Steffy$^3$}
  \ZTPPreprint
  \ZTPNumber{17-43}
  \ZTPMonth{July}
  \ZTPYear{2017}
  \ZTPInfo{
  $^1$ Zuse Institute Berlin, Department of Mathematical Optimization, Takustr.~7, 14195~Berlin, Germany, \nolinkurl{miltenberger@zib.de}\\
  $^2$ Department of Industrial and System Engineering, Lehigh University, Bethlehem, PA 18015 USA, \nolinkurl{ted@lehigh.edu}\\
  $^3$ Mathematics and Statistics, Oakland University, Rochester, MI 48309 USA, \nolinkurl{steffy@oakland.edu}  
  }
}
{
}

\begin{document}
\ifthenelse{\zibreport = 1}{}
{
\mainmatter
\titlerunning{Numerics in MILP}
\authorrunning{Miltenberger, Ralphs, Steffy}

\tocauthor{Matthias Miltenberger, Ted Ralphs, Dan Steffy}

\institute{
Zuse Institute Berlin, 14195 Berlin, Germany,
\email{miltenberger@zib.de},
\and
Department of Industrial and System Engineering, Lehigh University, Bethlehem,
PA 18015 USA,
\email{ted@lehigh.edu},
\and
Mathematics and Statistics, Oakland University, Rochester, MI 48309 USA,
\email{steffy@oakland.edu }
}
\author{Matthias Miltenberger\inst{1} \and Ted Ralphs\inst{2} \and Daniel 
E. Steffy\inst{3}}
}

\title{Exploring the Numerics of Branch-and-Cut for Mixed Integer Linear Optimization}
\date{July 20, 2017}

\ifthenelse{\zibreport = 1}{
\author{Matthias Miltenberger
\and Ted Ralphs
\and Daniel E. Steffy
}

\zibtitlepage
\thispagestyle{empty}
\setcounter{page}{0}
}{}

\maketitle

\begin{abstract}
We investigate how the numerical properties of the LP relaxations evolve
throughout the solution procedure in a solver employing the branch-and-cut
algorithm. The long-term goal of this work is to determine whether the effect
on the numerical conditioning of the LP relaxations resulting from the
branching and cutting operations can be effectively predicted
and whether such predictions can be used to make better algorithmic
choices. In a first step towards this goal, we discuss here the numerical
behavior of an existing solver in order to determine whether our 
intuitive understanding of this behavior is correct. 
\ifthenelse{\zibreport = 1}{}
{
}
\end{abstract}

\section{Introduction}
The branch-and-cut algorithm for mixed integer linear optimization problems
(MILPs) combines aspects of the branch-and-bound algorithm with the cutting
plane algorithm to strengthen the initial LP relaxation (see~\cite{Conforti2014} for a complete description of these operations and the definitions of these terms).
While branching increases the number of subproblems to be solved and should thus be avoided in principle, the addition of too many cutting planes often results in an LP relaxation with undesirable numerical properties.
Recent research into the viability of solving MILPs using a pure cutting plane approach has provided some insight into how and why this happens and has explored techniques to generate a sequence of valid inequalities whose addition to the LP relaxation is less likely to cause difficulties~\cite{Fischetti2011,Zanette2011}.

In general, branching and cutting must be used carefully in concert with each
other to maintain numerical stability. The effect of these operations on
numerics is not well understood, however, and is difficult to control
directly. There exists a number of approaches to effectively combine the
branching and cutting operations. In some solvers, cutting is only done at the
root node, while in others, cuts are added throughout the tree. As with any
numerical process, implementations of these solution algorithms use floating-point arithmetic and are subject to accumulation of
roundoff errors within the computations. Without appropriate handling of these
errors, the algorithms may return unreliable results, failing to behave or
terminate as expected.

Modern MILP solvers use a wide range of techniques to mitigate the
difficulties associated with numerical errors. For example, it is standard
practice
to discard or modify cuts whose coefficients differ significantly in
magnitude, since these inequalities are likely to degrade the conditioning of the LP
relaxation. This and other techniques help to ensure that the LP relaxation
will have better numerical properties and increases the computational stability
of the algorithm.

%

It is well understood that the addition of cutting planes has the potential to
negatively impact the numerical properties of the LP relaxation, even after
steps have been taken to improve their reliability. On the other hand,
branching may counteract this effect to some extent, leading to a more stable
algorithm overall.
In this paper we seek to carefully investigate the impact of both branching and cutting on the numerical properties of the LP subproblems solved in the branch-and-cut algorithm.
The purpose of this work is both to confirm existing folklore, namely that
branching improves condition and cutting degrades it, as well as to explore
the potential for directly controlling numerical properties through judicious
algorithmic choices. 

In Section~\ref{sec.condition} we discuss the choice and computation of the basis matrix condition numbers as a measure of numerical stability.
In Section~\ref{sec.computations} we describe computational results regarding how branching and cutting affect the condition numbers.
Section~\ref{sec.outlook} discusses some implications of our findings and ongoing work.

\section{Condition numbers}\label{sec.condition}

The \emph{condition number} of a numerical problem is a bound on the relative
change (in terms of a given norm) in the solution to a problem that can occur
as a result of a change in the input (see~\cite{BuergisserCucker2013} for formal definitions).
For example, the condition number of a matrix $A$ is $\kappa(A)=\|A\|_2\|A^{-1}\|_2$ and yields a bound on how much the solution to the linear system of equations $Ax=b$
might change, relative to a change in the right hand side vector $b$. For LPs,
a handful of different condition numbers have also been defined; a
comprehensive treatment of condition numbers for LPs, along with much more
general discussion regarding the concept of problem condition, is given in
\cite{BuergisserCucker2013}.

When LPs are solved by the simplex method, a sequence of basis matrices are
encountered (see~\cite{Conforti2014}), each corresponding to a square system of linear
equations. Although condition numbers can be defined for LPs themselves, it is 
the condition number of the basis matrices encountered during a simplex solve
(particularly the optimal basis) that is the most relevant measure of
numerical stability of the branch-and-cut algorithm. A primary reason for this
is that the solution to the LP relaxation is obtained by solving a system of
equations involving the basis matrix so that the condition number of this
matrix determines the multiplicative effect of numerical errors in the
computed cuts.

After applying cutting planes or branching at a node, the resulting modified
LP is re-solved.
In general, we
expect that the newly added cuts or branching inequalities will be binding 
at the new basic solution, which means that these additional constraints are
a factor in determining the conditioning of the basis. 
Thus, measuring the condition number of these linear systems and how they
change as a result of the added cuts or branching inequalities should give
some insight into the numerical behavior of the simplex algorithm and
ultimately the branch-and-cut algorithm. In this paper, we are looking for
overall trends (how much does the addition of cuts \emph{generally} degrade
the conditioning), so we consider these numbers in the aggregate and provide
some suggestions for visualizing this data.



Since we are are interested in an accurate picture, we use the 2-norm power
iteration method to determine condition numbers. This method provides an
accurate answer, though it is unlikely to be efficient enough for practical
use.
An excellent discussion on algorithms for condition number estimation is given in \cite[Chapter 15]{Higham2002}.

\section{Experiments}\label{sec.computations}

To study the effect of cuts on conditioning, we solved a subset of instances
from MIPLIB 3\cite{BixbyCeriaMcZealSavelsbergh98}, MIPLIB
2003\cite{AchterbergKochMartin06}, and MIPLIB 2010\cite{KochEtAl11} test
sets, collecting detailed statistics. The solver used was SCIP 4.0 with the LP
solver SoPlex 3.0~\cite{MaherFischerGallyetal.2017} (with slight modifications to allow access to the condition
number information).
We used a time limit of one hour and a node limit of 10,000.


To get a clearer picture, we deactivated many advanced features, such as
primal heuristics, domain propagation, and conflict analysis. Furthermore, we
only generated Gomory cuts and disabled all other cutting plane generators.
While SCIP only applies cutting planes at the root by default, we enabled cut
generation at all nodes in order to study how this affects conditioning. Note
that although cuts are generated throughout the tree, SCIP still uses a scoring
strategy to determine which inequalities should actually be added.


In what follows, we first study how the condition number of the basis matrix
evolves at the root node, where the initial LP relaxation is solved and initial
rounds of cuts are added,
and then study how the condition number of the basis matrices are affected by
branching and cutting as the algorithm progresses. 

\subsection{Root Node Analysis}

\begin{figure}
\vspace{-.2cm}
\centering
\includegraphics[width=.95\textwidth]{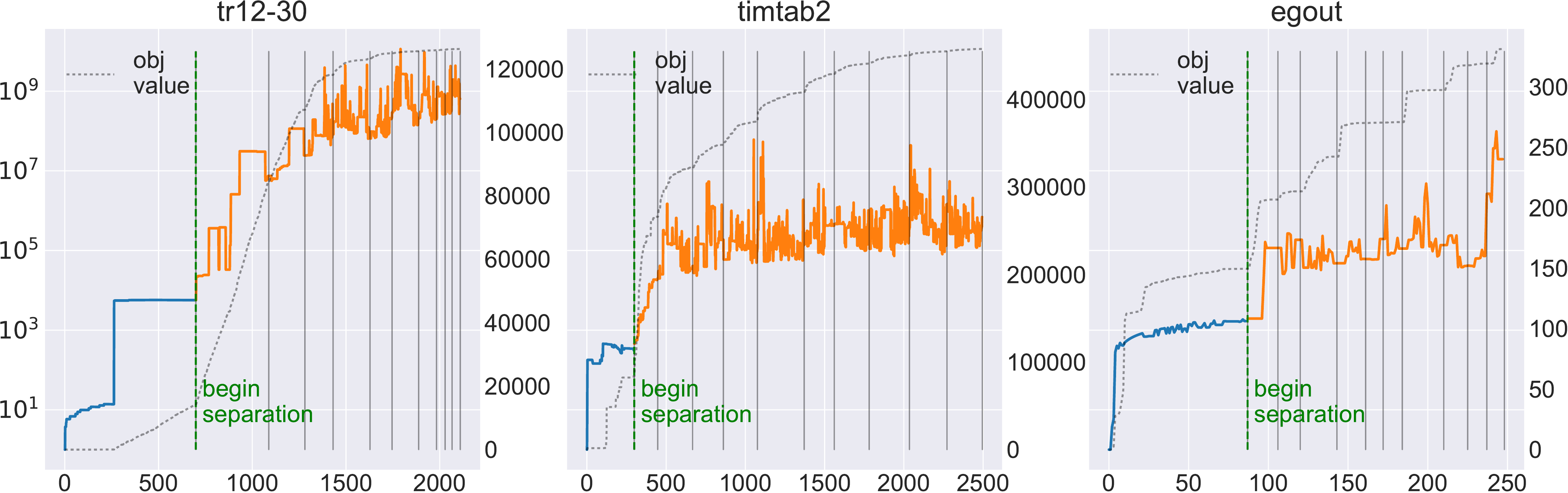}
\caption{Condition number development (vertical axis, in log scale) for every simplex iteration in the root node (horizontal axis) including re-optimizations after adding cutting planes in multiple rounds (vertical lines). A plot of objective values at each iteration is overlaid as a dashed gray line with the scales given to the right of each plot.}
\label{fig:simplex-condition-numbers-root}
\vspace{-.2cm}
\end{figure}

In general, we expect the condition number of the basis matrix to degrade as a
result of operations performed in the root node and our initial computations
are aimed at confirming this. Fig.~\ref{fig:simplex-condition-numbers-root}
shows the condition number of each basis matrix encountered during each
iteration of the solution of the initial LP relaxation and during each
iteration of the re-solve occurring after adding each round of cuts for
selected instances from our test set.
One can observe that during the early iterations---especially of the initial
relaxation in the root node---the condition numbers of the basis matrices
grow quickly. This is expected, as more structural variables are pivoted into the
basis, while slack variables are pivoted out. Since the initial basis is always
the identity matrix, which has condition number 1, the conditioning can only
degrade at first. After the initial optimization, the MILP solver tries to
generate Gomory cuts.
\begin{wrapfigure}{r}{.4\textwidth}
\vspace{-.4cm}
\centering
\includegraphics[width=.9\linewidth, height=.9\linewidth]{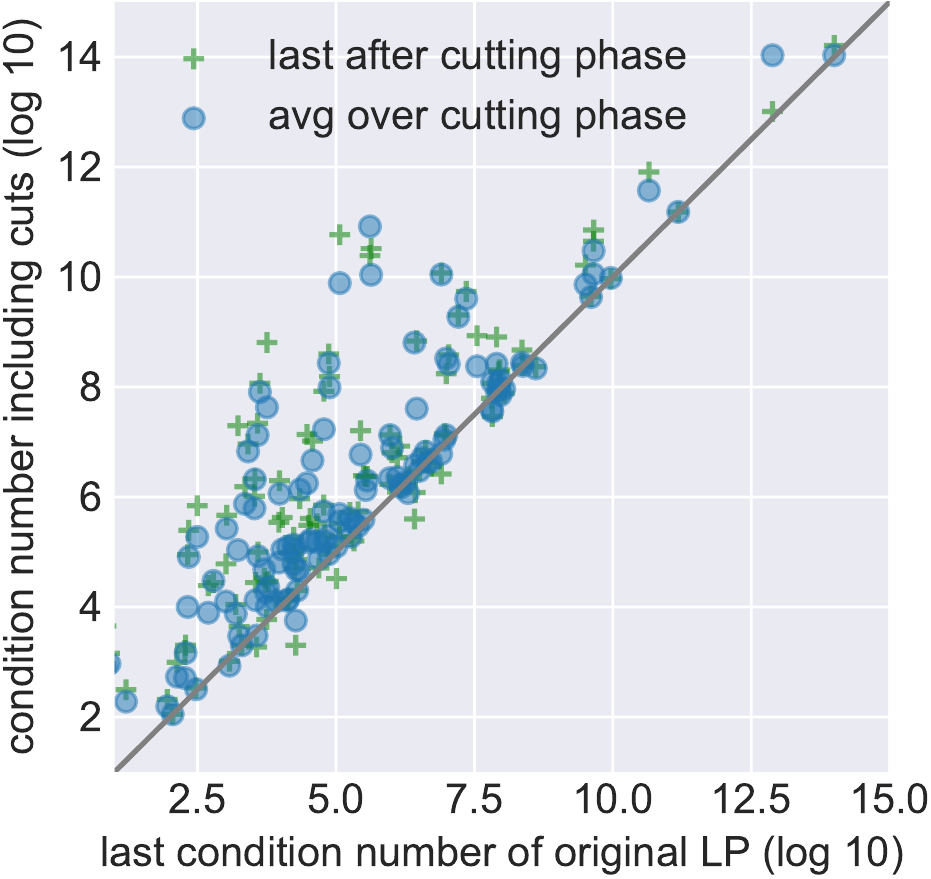}
\caption{Root node: Comparison of condition numbers of the original LP and including cutting planes.}
\vspace{-0.5cm}
\label{fig:orig-vs-cuts-root}
\end{wrapfigure}
This computation involves the basis matrix itself, so an ill-conditioned basis
matrix can prevent precise calculation of the coefficients of the new
constraint. Moreover, adding these new rows to the LP often deteriorates its
condition number even further as can be seen in
Fig.~\ref{fig:simplex-condition-numbers-root}. This sample of instances
clearly shows the expected behavior.

Fig.~\ref{fig:orig-vs-cuts-root} is a visualization of the difference
between the condition number of the optimal basis of the original LP and two
other numbers: (1) the average over all bases encountered during the cutting
procedure and (2) the condition number of the final optimal basis. While for
some instances there is a slight improvement after adding cuts, in most cases
addition of cuts leads to an increased condition number, as expected.

\subsection{Tree Analysis}

One way in which the addition of cuts can cause basis matrices to become
poorly conditioned is if the associated hyperplanes are nearly parallel;
addition of many such cuts may lead to a \emph{tailing off} of the cutting plane algorithm as many similar cuts are generated and
the process stalls. Although branching also involves imposing a special kind
of ``cut'' to the resulting subproblems, these branching constraints have a
simple form (the coefficient vector is a unit vector), which makes them quite
attractive from a numerical point of view. In particular, they are mutually
orthogonal and unlikely to degrade the conditioning much in general. As such,
we may be tempted to hope that the addition of this special kind of
inequality may even improve conditioning.



\begin{wrapfigure}{l}{.3\textwidth}
\vspace{-.5cm}
\centering
\includegraphics[width=2.5cm]{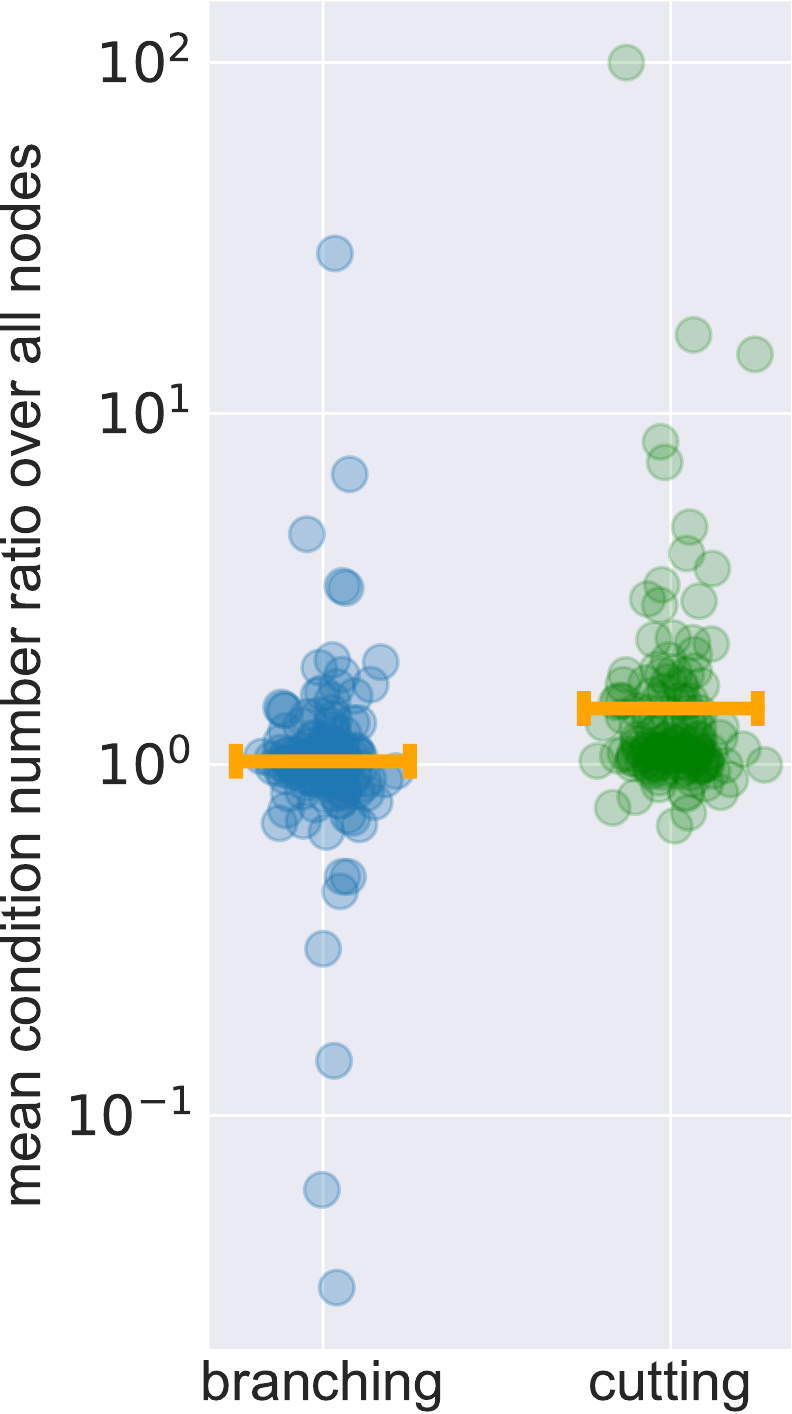}
\caption{Effects of branching and cutting}
\label{fig:branching-vs-cutting-change}
\vspace{-.3cm}
\end{wrapfigure}
Despite the apparent plausibility of this hypothesis, our experiments do not
fully support it, though they do show a significant difference between the
effect of branching versus cutting, as expected. In
Fig.~\ref{fig:branching-vs-cutting-change}, we show how branching and
cutting impact the numerical stability. The left plot shows the average
relative change in the condition number as a result of the addition of the
branching constraints.
Similarly, the right plot shows the average relative change in conditioning resulting
from the addition of cuts.
In each case, we took the difference between the condition
numbers of the optimal basis matrices before and after either branching or
cutting. Each dot then represents the average across all nodes for a given
instance.
The bar represents the mean over all instances.
While branching does not seem to have a significant effect on average, adding
cutting planes clearly leads to an increase in the condition number.
Thus, despite the observation that branching does not appear to degrade the condition number in the same way as cut generation, it does not appear to help it either.


\begin{figure}
\vspace{-.5cm}
\centering
\includegraphics[width=.9\linewidth]{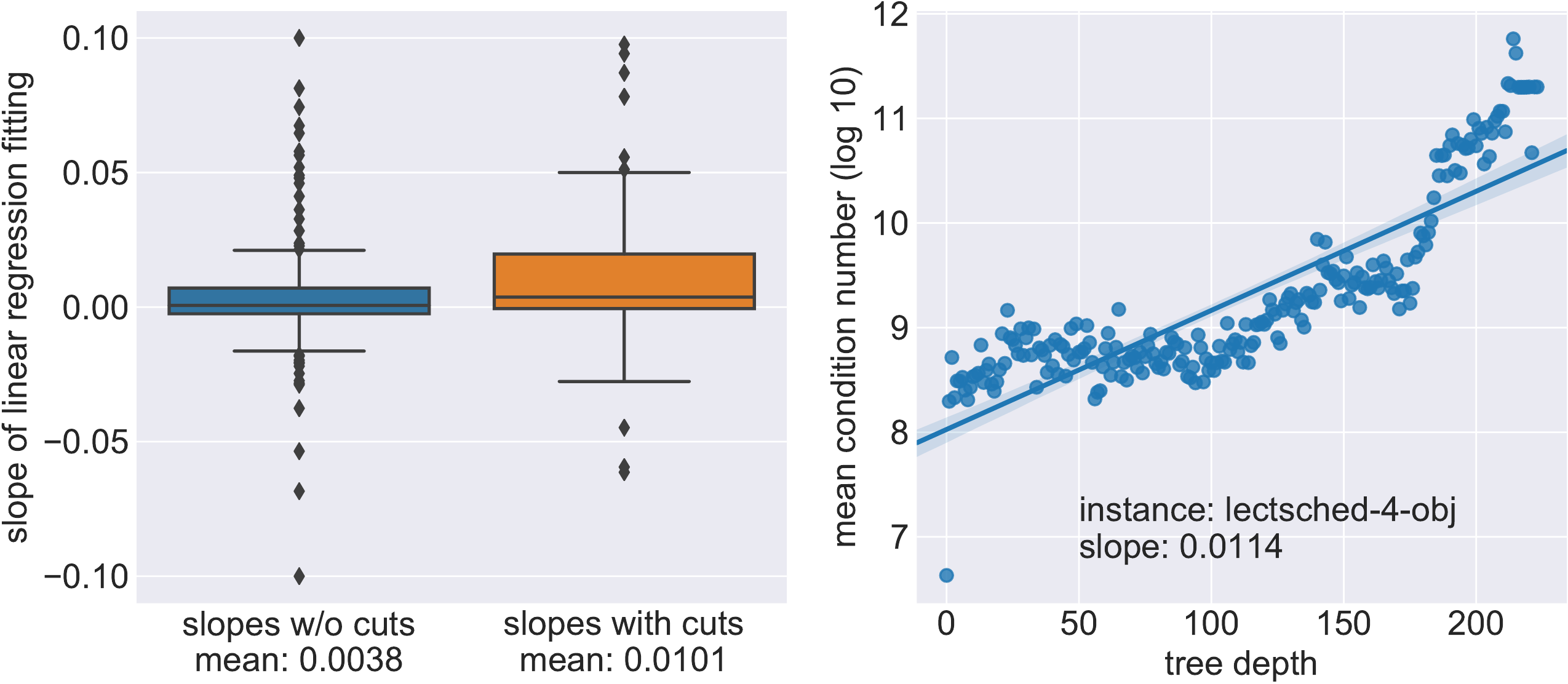}
\caption{Condition number development in the tree. Left: Distribution of linear regression slopes of all instances in the test set. Right: Single instance example.}
\label{fig:tree_regression}
\vspace{-.5cm}
\end{figure}

In Fig.~\ref{fig:tree_regression} we visualize how condition
numbers degrade generally as a function of the depth of a given node. The idea
is to determine whether conditioning generally degrades consistently as the
tree gets deeper. The right figure plots the average condition number across
all nodes at a given depth, along with a regression line showing the average
degradation in the log of the condition number per level in the tree for a
single instance. The left figure shows the distribution of slopes of this same
linear regression across all instances both with cuts and with a pure branch-and-bound.

It appears that in general, the condition number often has a strong positive correlation with the tree depth if cuts are added throughout the solving process.
When cutting is disabled this effect is much less strong.
One has to be aware that the behavior of a single instance might be much different from what the trend predicts.

\section{Outlook}\label{sec.outlook}

In this paper, we presented a preliminary exploration of the numerical
behavior of SCIP, a state-of-the-art MILP solver. In the future, we hope to do
similar explorations with other solvers to determine what the overall behavior
is and where additional control of the numerical stability might have an
impact. The eventual goal is to determine whether it is possible to more
directly estimate the impact of certain algorithmic choices on numerical
behavior and whether this could lead to improved control mechanisms.

\subsubsection*{Acknowledgments}
The work for this article has been partly conducted within the \emph{Research Campus Modal} funded by the German Federal Ministry of Education and
Research (fund number 05M14ZAM). The support of Lehigh University is also
gratefully acknowledged.

\bibliographystyle{abbrv}
\bibliography{numerics_in_mip}

\begin{thebibliography}{1}

\bibitem{AchterbergKochMartin06}
T.~Achterberg, T.~Koch, and A.~Martin.
\newblock {MIPLIB} 2003.
\newblock {\em Operations Research Letters}, 34(4):1--12, 2006.

\bibitem{BixbyCeriaMcZealSavelsbergh98}
R.~E. Bixby, S.~Ceria, C.~M. McZeal, and M.~W.~P. Savelsbergh.
\newblock An updated mixed integer programming library: {MIPLIB} 3.0.
\newblock {\em Optima}, (58):12--15, June 1998.

\bibitem{BuergisserCucker2013}
P.~B{\"{u}}rgisser and F.~Cucker.
\newblock {\em Condition - The Geometry of Numerical Algorithms}, volume 349 of
  {\em Grundlehren der math. Wissenschaften}.
\newblock Springer, 2013.

\bibitem{Conforti2014}
M.~Conforti, G.~Cornu\'ejols, and G.~Zambelli.
\newblock {\em Integer Programming}.
\newblock Springer, 2014.

\bibitem{Fischetti2011}
M.~Fischetti and D.~Salvagnin.
\newblock A relax-and-cut framework for gomory mixed-integer cuts.
\newblock {\em Math Prog Comp}, 3(2):79--102, 2011.

\bibitem{Higham2002}
N.~J. Higham.
\newblock {\em Accuracy and Stability of Numerical Algorithms}.
\newblock 2002.

\bibitem{KochEtAl11}
T.~Koch, T.~Achterberg, E.~Andersen, O.~Bastert, T.~Berthold, R.~E. Bixby,
  E.~Danna, G.~Gamrath, A.~M. Gleixner, S.~Heinz, A.~Lodi, H.~Mittelmann,
  T.~Ralphs, D.~Salvagnin, D.~E. Steffy, and K.~Wolter.
\newblock {MIPLIB} 2010.
\newblock {\em Math Prog Comp}, 3(2):103--163, 2011.

\bibitem{MaherFischerGallyetal.2017}
S.~J. Maher, T.~Fischer, T.~Gally, G.~Gamrath, A.~Gleixner, R.~L. Gottwald,
  G.~Hendel, T.~Koch, M.~E. L{\"u}bbecke, M.~Miltenberger, B.~M{\"u}ller, M.~E.
  Pfetsch, C.~Puchert, D.~Rehfeldt, S.~Schenker, R.~Schwarz, F.~Serrano,
  Y.~Shinano, D.~Weninger, J.~T. Witt, and J.~Witzig.
\newblock {The SCIP Optimization Suite 4.0}.
\newblock Technical Report 17-12, ZIB, 2017.

\bibitem{Zanette2011}
A.~Zanette, M.~Fischetti, and E.~Balas.
\newblock Lexicography and degeneracy: can a pure cutting plane algorithm work?
\newblock {\em Math Prog}, 130(1):153--176, 2011.

\end{thebibliography}

\end{document}